\documentclass[sort&compress]{elsarticle}
\usepackage{graphicx}
\usepackage{amsmath}
\usepackage{amssymb}
\usepackage{color}
\usepackage{amsthm}
\usepackage[english]{babel}
\usepackage{setspace}

\usepackage{mathtools}
\usepackage{xparse}

\newtheorem{theorem}{Theorem}

\newtheorem{lemma}[theorem]{Lemma}
\newtheorem{definition}[theorem]{Definition}
\newtheorem{proposition}[theorem]{Proposition}

\newcommand{\cB}{\mathcal{B}}

\newcommand{\bai}{\bar{i}}
\newcommand{\baj}{\bar{j}}
\newcommand{\bak}{\bar{k}}

\renewcommand{\le}{\leqslant}
\renewcommand{\leq}{\leqslant}
\renewcommand{\ge}{\geqslant}
\renewcommand{\geq}{\geqslant}

\newcommand{\Z}{\mathbb{Z}}

\newcommand{\kp}{{k_+}}
\newcommand{\km}{{k_-}}
\newcommand{\ve}{\mathbf{e}}

\newcommand{\vv}{\mathbf{v}}

\newcommand{\vx}{\mathbf{x}}

\newcommand{\eqdef}{\triangleq}
\newcommand{\splt}{\diamond}

\DeclarePairedDelimiterX{\set}[1]{\{}{\}}{\setargs{#1}}
\NewDocumentCommand{\setargs}{>{\SplitArgument{1}{;}}m}
{\setargsaux#1}
\NewDocumentCommand{\setargsaux}{mm}
{\IfNoValueTF{#2}{#1} {#1\, : \,\mathopen{}#2}}

\DeclarePairedDelimiter\abs{\lvert}{\rvert}

\DeclarePairedDelimiter\floor{\lfloor}{\rfloor}
\DeclarePairedDelimiter\parenv{\lparen}{\rparen}

\begin{document}

\title{On Lattice Tilings of Asymmetric Limited-Magnitude Balls $\cB(n,2,m,m-1)$\tnoteref{t1}}
\tnotetext[t1]{This work was supported in part by the National Natural Science Foundation of China under Grant 12371523.}
\author[1]{Zhihao Guan}
\ead{guanzh@stu.xjtu.edu.cn}

\author[1]{Hengjia Wei\corref{cor1}}
\ead{hjwei05@gmail.com}

\author[2]{Ziqing Xiang}
\ead{xiangzq@sustech.edu.cn}

\cortext[cor1]{Corresponding author}
\address[1]{ School of Mathematics and Statistics, Xi'an Jiaotong University, Xi'an 710049, China}
\address[2]{ Department of Mathematics and National Center for Applied Mathematics  Shenzhen,  Southern University of Science and Technology, Shenzhen 518055, China}

\date{}

\begin{abstract}Limited-magnitude errors modify a transmitted integer vector in at most $t$ entries, where each entry can increase by at most $\kp$ or decrease by at most $\km$. This channel model is particularly relevant to applications such as flash memories and DNA storage. A perfect code for this channel is equivalent to a tiling of $\Z^n$ by asymmetric limited-magnitude balls 
$\cB(n,t,\kp,\km)$. In this paper, we focus on the case where $t=2$ and $\km=\kp-1$, and we derive necessary conditions on $m$ and $n$ for the existence of a lattice tiling of $\cB(n,2,m,m-1)$.  Specifically, we prove that if such a tiling exists, then either $4\leq m \leq 512$ and $n<7.23m+4$, or $m>512$ and $n<4m$. In particular, for $m=2$ and $m=3$, we show that no lattice tiling of $\cB(n,2,2,1)$ or $\cB(n,2,3,2)$ exists for any $n\geq 3$.
\end{abstract}

\maketitle

\section{Introduction}
In certain applications, such as flash memories \cite{cassuto2010codes,schwartz2011quasi} and DNA based data storage systems \cite{Jainetal2020ISIT,LeeKalGoeBolChur19,WeiSch21}, information is encoded as an integer vector $x\in \Z^n$.  These applications are typically affected by noises that are characterized by limited magnitude errors. Specifically, at most $t$ entries of the transmitted vector can be altered, with each entry changing by at most $\kp$ in the positive direction or by at most $\km$ in the negative direction. For integers $n \geq t \geq 1$ and $\kp \geq \km \geq 0$, we define the \emph{$(n,t,\kp,\km) $-error-ball} as
$$
\cB(n,t,\kp,\km) \triangleq \left\{ \vx = (x_1, x_2, \dots, x_n)\in \Z^n:-\km \leq x_i \leq \kp \text{ and } \text{wt}(\vx) \leq t \right\},
$$
where wt$(\vx)$ denotes the Hamming weight of $\vx$. It is readily checked that an error-correcting code in this setting is equivalent to a packing of $\Z^n$ by the error-balls $\cB(n,t,\kp,\km)$ and a perfect code is equivalent to a tiling of $\Z^n$ by $\cB(n,t,\kp,\km)$. 

Previous work on tiling these shapes has focused primarily on the case of $t=1$. The cross, $\cB(n,1,k,k)$, and the semi-cross, $\cB(n,1,k,0)$, have been extensively studied in \cite{hamaker1984combinatorial,hickerson1986abelian,klove2011some,stein1984packings,stein1994algebra}. This research was later extended to quasi-crosses, $\cB(n,1,\kp,\km)$, in \cite{schwartz2011quasi}, sparking a surge of activity on the subject \cite{schwartz2014non,yari2013some,ye2019some,zhang2016new,zhang2017nonexistence,zhang2017splitter}. Additionally, crosses with arms of length half have been studied in \cite{BuzEtz13}.

However, for the case where $t\geq 2$, there is a limited body of work. Tilings of $\Z^n$ by $\cB(n,n-1,k,0)$ have been studied in \cite{buzaglo2012tilings,klove2011some,stein1990notched}. In \cite{wei2021lattice,WEI2022103450}, the authors initiated the study of the tiling and packing problems for the general cases where\footnote{When $t=n$, $\cB(n,n,\kp,\km)$ is a hypercube, which can tile  $\Z^n$ by $((\kp+\km+1)\Z)^n$.} $2\leq t\leq n-1$. A lattice tiling construction based on perfect codes in the Hamming metric was given, and several non-existence results were proven. 
In particular, in \cite{wei2021lattice}, it is shown that $\cB(n,2,1,0)$ lattice-tiles $\Z^n$ if and only if $n\in\{3,5\}$, and $\cB(n,2,2,0)$ lattice-tiles $\Z^n$ if and only if $n=11$. Subsequently, in \cite{zhang2023lattice}, it is shown that $\cB(n,2,1,1)$ lattice-tiles $\Z^n$ if and only if $n=11$.

In this paper, we consider the case where $t=2$ and $\kp=\km+1$. Denote $m\eqdef \kp=\km+1\ge 2$. We study the necessary conditions on $n$ such that the ball $\cB(n,2,m,m-1)$ lattice-tiles $\Z^n$. The main results of this paper are as follows.

\begin{theorem}\label{thm:main}
Let $n\geq 3$ and $m\geq 2$. If  $\cB(n,2,m,m-1)$ lattice-tiles $\mathbb{Z}^n$, then either $4\leq m \leq 512$ and $n<7.23m+4$, or $m>512$ and $n<4m$.
\end{theorem}
It follows from the theorem above that for $n\geq 3$, there are no lattice tilings of $\Z^n$ by $\cB(n,2,2,1)$ or $\cB(n,2,3,2)$.

The proof of Theorem~\ref{thm:main} proceeds as follows. It is known that a lattice tiling of 
$\cB$ is equivalent to the existence of an Abelian group  $G$ satisfying certain combinatorial conditions. In this paper, we analyze the quantity 
$m_2(G)$, defined as the number of elements of order $2$ in $G$ plus 1. We derive both lower and upper bounds on $m_2(G)$, which are then used to establish Theorem \ref{thm:main}.


This paper is organized as follows. Section~\ref{Sec:Pre} introduces the notations and fundamental concepts used throughout the paper, including lattice tiling, group splitting, and their equivalence. Section~\ref{Sec:Lower} and Section~\ref{Sec:Upper} present several lower bounds and upper bounds for $m_2(G)$, respectively. Section~\ref{Sec:main} completes the proof of Theorem~\ref{thm:main}.

\section{Preliminaries}\label{Sec:Pre}
Let $n,t$ be integers such that $n\ge t\ge 1$. Let $\kp$, $\km$ be integers such that $\kp \geq \km \geq 0$. For integers $a\le b$, we denote $[a,b]\triangleq \{a,a+1,\ldots,b\}$ and $[a,b]^* \triangleq  [a,b] \setminus \{0\}$.

A \textit{lattice} $\Lambda\subseteq \Z^n$ is an additive subgroup of $\Z^n$. A lattice $\Lambda$ may be represented by a matrix $\mathcal{G}(\Lambda)\in \Z^{n\times n}$, the span of whose rows (with integer coefficients) is $\Lambda$. A \textit{fundamental region} of $\Lambda$ is defined as
$$\left\{\sum_{i=1}^n c_i \vv_i:c_i\in \mathbb{R},0\leq c_i<1\right\},$$
where $\vv_i$ is the $i$th row of $\mathcal{G}(\Lambda)$. The volume of the fundamental region is $|\det \mathcal{G}(\Lambda)|$, and is independent of the choice of $\mathcal{G}(\Lambda)$.

We say $\cB\subseteq \Z^n$ \textit{packs} $\Z^n$ by $\Lambda\subseteq\Z^n$, if the translates of $\cB$ by elements from $\Lambda$ do not intersect, i.e., for all $\vv$, $\vv'\in \Lambda$, $\vv\neq \vv'$,
$$(\vv+\cB)\cap (\vv'+\cB)=\varnothing.$$
We say $\cB$ \textit{covers} $\Z^n$ by $\Lambda$ if
$$\bigcup_{\vv\in\Lambda}(\vv+\cB)=\Z^n.$$
If $\cB$ both packs and covers $\Z^n$ by $\Lambda$, then we say $\cB$ \textit{tiles} $\Z^n$ by $\Lambda$.

An equivalence between the lattice packing of $\cB(n,1,\kp,\km)$ and group splitting is described in \cite{hickerson1986abelian,stein1984packings}. In \cite{schwartz2011quasi}, the technique is used to solve the quasi-cross lattice tiling problem. For $t\geq 2$, there is an extended definition of group splitting. 

\begin{definition}
Let $G$ be a finite Abelian group, where $+$ denotes the group operation. For $m\in\Z$ and $g \in G$, let $mg$ denote $g +g +\cdots+g$ (with $m$ copies of $g$) when $m>0$, which is extended naturally to $m \leq 0$. Let $M \subseteq \Z\setminus \{0\}$ be a finite set, and $S=\{s_1,s_2,\ldots,s_n\}\subseteq G$. 
\begin{itemize}
    \item If the elements $\ve\cdot(s_1,\ldots,s_n)$, where $\ve \in (M \cup \{0\})^n$ and $1\leq \text{wt}(\ve)\leq t$, are all distinct and non-zero in $G$, we say the set $M$ \textit{partially $t$-splits} $G$ with splitter set $S$, denoted
    $$G\geq M\splt_t S.$$
    \item If for every $g\in G$ there exists a vector $\ve \in (M \cup \{0\})^n$ with wt$(\ve)\leq t$, such that $g=\ve\cdot(s_1,\ldots,s_n)$, we say the set $M$ \textit{completely $t$-splits} $G$ with splitter set $S$, denoted 
    $$G\leq M\splt_t S.$$
    \item If $G\geq M\splt_t S$ and $G\leq M\splt_t S$ we say $M$ \textit{$t$-splits} $G$ with splitter set $S$, and write
    $$G=M\splt_t S.$$
\end{itemize}
\end{definition}

Intuitively, $G=M\splt_t S$ means that the non-trivial linear combinations of elements from $S$, with at most $t$ non-zero coefficients from $M$, are distinct and give all the non-zero elements of $G$ exactly once. 

The following two theorems show the equivalence of $t$-splittings and lattice tilings, summarizing Lemma 3, Lemma 4, and Corollary 1 in \cite{buzaglo2012tilings}. They generalize the treatment for $t=1$ in previous work.

\begin{theorem}[Lemma 4 and Corollary 1 in \cite{buzaglo2012tilings}]\label{thm:gs2lt}
Let $G$ be a finite Abelian group, $M\triangleq [-\km,\kp]^*$, and
$S=\{s_1,\ldots,s_n\} \subseteq G$, such that $G =M\splt_t S$. Define $\phi: \Z^n \to G$ as $\phi(\vx) \triangleq \vx\cdot(s_1,\ldots,s_n)$ and let
$\Lambda\triangleq\ker\phi$ be a lattice. Then $\cB(n,t,\kp,\km)$ tiles $\Z^n$ by $\Lambda$.
\end{theorem}

\begin{theorem}[Lemma 3 and Corollary 1 in \cite{buzaglo2012tilings}]\label{thm:lt2gs}
Let $\Lambda\subseteq \Z^n$ be a lattice, and assume $\cB(n,t,\kp,\km)$ tiles $\Z^n$ by $\Lambda$. Then there exists a finite Abelian group $G$ and $S=\{s_1,\ldots,s_n\} \subseteq G$ such that $G =M\splt_t S$, where $M\triangleq [-\km,\kp]^*$.
\end{theorem}

In this paper, we study the lattice tiling problem for $\cB(n,2,m,m-1)$.  
By Theorem~\ref{thm:gs2lt} and Theorem~\ref{thm:lt2gs}, this is equivalent to examining the $2$-splitting problem with $M\eqdef [-m+1,m]^*$. The following lemma counts the number of pairs in $M$ that sum to $x$, which will be used later.
\begin{lemma}\label{lemma:tau} For each $x \in [-2(m-1),2m]$, let   
    \begin{align*}
    \tau(x) \eqdef \left|\set{(a,b)\in M^2; a+b=x}\right|.
    \end{align*}
    Then the following results hold:

    1) $\tau(x)=
    \begin{cases}
        2m+1-x  & m+1\leq x \leq 2m, \\
        2m-1-x & 1\leq x \leq m, \\
        2m-2 & x=0, \\
        2m-3+x & -m+1 \leq x \leq -1, \\
        2m-1+x & -2m+2\leq x \leq -m;
    \end{cases}$
    
    2) $\sum\limits_{x=m+1}^{2m} \tau(x)=\frac{m^2+m}{2}$;

    3) $\sum\limits_{x=-2m+2}^{-m} \tau(x)=\frac{m^2-m}{2}$;

    4) $\sum\limits_{\substack{x=-2m+3\\ x \text{odd}}}^{2m-1} \tau(x)=2m^2-2m$.
\end{lemma}
\begin{proof}
    Part 1) can be obtained directly by counting. Parts 2), 3) and 4) are derived from the calculations in part 1).
\end{proof}

To deal with the case of $m=2$, we require the following result regarding the Ramanujan-Nagell equation.

\begin{lemma}[Nagell, \cite{nagell1961diophantine}]\label{lemma:Nagell}
    The only solutions $(x,n)$ in natural numbers of the equation 
    $$2^n-7=x^2$$
    are $(1,3),(3,4),(5,5),(11,7),(181,15)$.
\end{lemma}





\section{Lower bounds on $m_2(G)$}\label{Sec:Lower}

Let $G$ be an Abelian group.
For a non-zero element $g\in G$,  define
\[ord(g)   \triangleq \min\{n\in \mathbb{Z}: n>0, ng=0\}.\]
For an integer $l\geq 2$,  define
\[m_l(G)  \triangleq |\{g\in G: ord(g)=l\}|.\]
We are going to derive two lower bounds on $m_2(G)$, stated in Lemma~\ref{lm:lowerbound1} and Lemma~\ref{lm:lowerbound2}, using a counting method adapted from \cite{WEI2022103450}. The first bound is independent of $m$, while the second incorporates $m$  and is particularly useful for small values of $m$. Subsection~\ref{subsec:counting} presents the counting framework, and Subsection~\ref{subsec:lbndpf} provides the proofs of Lemma~\ref{lm:lowerbound1} and Lemma~\ref{lm:lowerbound2}.

\begin{lemma}\label{lm:lowerbound1}
    Let $n\geq 4m$.  Suppose that $G=[m+1,m]^*\splt_2 S$. Then
    $$m_2(G)\geq \frac{n}{2}+1.$$
\end{lemma}

\begin{lemma}\label{lm:lowerbound2}
    Let $n\geq m+4$. Suppose that $G=[m+1,m]^*\splt_2 S$. Denote $q'\triangleq \max\{d\mid |G|:\gcd(2m,d)=1\}$.
    \begin{enumerate}
        \item If $m_{2m}(G)>0$, then
            \begin{equation*}
            \begin{aligned}
            (2m-3)m_2(G)> & \ (2m-2)n^2-(4m^2-3)n -\left(\frac{n}{\sqrt{q'}}+1\right)^2\\
            & \ -(2mn-3n-2m^2-2m+6)\left(\frac{n}{\sqrt{q'}}+1 \right).\\
            \end{aligned}
            \end{equation*}
        \item If $m_{2m}(G)=0$, then
            $$(2m-3)m_2(G)\geq \left(2m-\frac{5}{2}\right)n^2-\left(4m^2-\frac{5}{2}\right)n.$$
    \end{enumerate}
\end{lemma}

\subsection{Counting argument}\label{subsec:counting}
Let $G$ be a finite Abelian group and assume that $G = [-m+1,m]^* \diamond_2 S$, for some $S=\{s_1,s_2,\ldots,s_n\}\subseteq G$. For $1\leq k\leq m-1$ and $1\leq i \leq n$, denote 
\[s_{i+kn}\triangleq (k+1)s_i  \ \textup{ and }  s_{i-kn}\triangleq -k s_i.\]
Furthermore, denote $s_{\infty}\triangleq 0$, and 
assume that $\infty \equiv \infty \pmod n$, $\infty \not\equiv i \pmod n$ and $i\not\equiv \infty \pmod n$ for all $i \in [1-n(m-1),nm]$. Let
$$\Delta \eqdef \set{(s_i,s_j);i,j\in [1-n(m-1),nm]\cup \set{\infty}, i\not\equiv j \pmod n }.$$
Then 
\begin{align}
\abs{\Delta} & =((2m-1)n+1)(2m-1)n-n(2m-1)(2m-2) \notag \\
& =(2m-1)^2 n^2-(4m^2-8m+3)n.\label{Delta}
\end{align}

We are going to estimate the number of distinct elements in the set $\set{s_i-s_j; (i,j)\in \Delta}$. 
This can be reduced to  estimating the number of the equations
\begin{equation} \label{countingEquation}
    s_i-s_j=s_k-s_l
\end{equation}
where $(s_i,s_j)$, $(s_k,s_l)\in \Delta$ and $(i,j)\neq (k,l)$. Note that Eq.~\eqref{countingEquation} implies
$$s_i+s_l=s_k+s_j.$$
Since $G\ge [-m+1,m]^* \diamond_2 S$, either $i\equiv l \pmod n$ or $k\equiv j \pmod n$. By exchanging the two sides of the equations, we assume that $i\equiv l \pmod n$ always holds. Let $M_1$ be the number of equations 
$s_i-s_j=s_k-s_l,$
where $(s_i,s_j)$, $(s_k,s_l)\in \Delta$, $i\equiv l \pmod n$ and $k\not\equiv j \pmod n$. 
Let $M_2$ be the number of equations
$s_i-s_j=s_k-s_l$,
where $(s_i,s_j)$, $(s_k,s_l)\in \Delta$, $i\equiv l \pmod n$ and $k\equiv j \pmod n$. 

We start with the set $\set{s_i-s_j; (i,j)\in \Delta}$. For each $\bai \in [1,n]\cup \set{\infty}$, if there are pairs $(i,j)$ and $(k,l)$ such that $s_i-s_j=s_k-s_l$, $i\equiv l \equiv \bai \pmod{n}$ and $j\not \equiv k \pmod{n}$, then we remove the element $(s_k,s_l)$ from $\Delta$. 
For any two $\bai,\baj \in [1,n]$ with $\bai < \baj$,  if there are pairs $(i,j)$ and $(k,l)$ such that $s_i-s_j=s_k-s_l$, $i\equiv l \equiv \bai \pmod{n}$ and $j \equiv k \equiv \baj \pmod{n}$, then we remove the element $(s_k,s_l)$ from $\Delta$. Finally, for the equations $s_{\infty
}-s_j = s_k-s_{\infty}$ with $j\equiv k \pmod{n}$, we remove the element $(s_k,s_{\infty})$ from $\Delta$. Denote the remaining set as $\Delta'$. Then $|\Delta'|\geq |\Delta|-M_1-M_2$. On the other hand, elements in the set $\set{s_i-s_j; (i,j)\in \Delta'}$ are non-zero and pairwise different. Then $|\Delta'|\leq |G|-1$. It follows that 
\begin{equation}\label{Inequality:Total}
    \begin{aligned}
        &\frac{(2m-1)^2}{2}n^2-\frac{4m^2-8m+3}{2}n = |G|-1 \ge |\Delta'| \\
        \ge & (2m-1)^2n^2-(4m^2-8m+3)n-M_1-M_2.
    \end{aligned}
    \end{equation}

In the following, we estimate the values of $M_1$ and $M_2$.

\begin{lemma}\label{lm:M1} Assume that $G \ge [-m+1,m]^* \diamond_2 S$.  Let $M_1$ be the number of equations 
$$s_i-s_j=s_k-s_l,$$
where $(s_i,s_j)$, $(s_k,s_l)\in \Delta$, $i\equiv l \pmod n$ and $k\not\equiv j \pmod n$. Then
    \begin{equation} \label{Maximum:Notequiv}
        M_1\leq 2m^2 n. 
    \end{equation}
\end{lemma}
\begin{proof}
    If $i=l=\infty$, then $s_j+s_k=0$. This implies that there are $a,b\in [-m+1,m]$ and distinct $\baj,\bak\in [1,n]$  such that $as_{\baj}+bs_{\bak}=0$, which contradicts the assumption that $G \ge [-m+1,m]^* \diamond_2 S$. 
    
    Now, let $\bar{i}$ be the unique integer of $[1,n]$ such that $i\equiv l\equiv \bar{i}$.  Since $s_i,s_l\in\{-(m-1)s_{\bar{i}},\ldots,-s_{\bar{i}},s_{\bar{i}},\ldots,m s_{\bar{i}}\}$, then $s_j+s_k=s_i+s_l\in \{-2(m-1)s_{\bar{i}},\ldots,2m s_{\bar{i}}\}$. Note that  $j\not \equiv \bar{i}$, $k\not \equiv \bar{i}$ and $k\not\equiv j \pmod n$. By the assumption that  $G \ge [-m+1,m]^* \diamond_2 S$, we have that  $s_j+s_k\notin \{-(m-1)s_{\bar{i}},\ldots,ms_{\bar{i}}\}$. For each $\bar{i}\in [1,n]$ and $x\in[-2(m-1),-m]\cup [m+1,2m]$, we claim that there is at most one unordered pair $\{j,k\}$ with $j\not \equiv k \pmod{n}$ such that $ s_j+s_k =x s_{\bar{i}}$. Otherwise, if we have another pair $\{j',k'\}$ satisfying the conditions, then $s_j+s_k =x s_{\bar{i}}=s_{j'}+s_{k'}$, a contradiction. Noting that a pair $\{j,k\}$ generates $2\tau(x)$ such equations, it follows that 
        
    $$M_1\leq n \left( \sum_{x=m+1}^{2m} 2\tau(x)+\sum_{x=-2m+2}^{-m} 2\tau(x)\right)=2 m^2 n,$$
    where the equation holds due to  Lemma \ref{lemma:tau}.
\end{proof}

Let 
    \[I\eqdef \{i\in [1,n]: 2m s_i=0 \}\] and denote 
    \[C\eqdef \abs{I}.\]
   Noting the $ms_i$'s, where $i \in I$, are distinct elements of order $2$,  we have $C\leq m_2(G)$. Meanwhile, for each $i\in I$, $ord(s_i)| 2m$, but $ord(s_i)> m$. Hence, $ord(s_i)=2m$. It follows that $C\leq m_{2m}(G)$. 

\begin{lemma}\label{lm:M2}
Assume that $G \ge [-m+1,m]^* \diamond_2 S$.  Let $M_2$ be the number of equations $$s_i-s_j=s_k-s_l,$$ where $(s_i,s_j)$, $(s_k,s_l)\in \Delta$, $i\equiv l \pmod n$ and $k\equiv j \pmod n$. 
\begin{enumerate}
\item If  $m_{2m}(G)>0$, then
    \begin{equation*}
    \begin{aligned}
       M_2 &\leq  \left(2m^2-4m+\frac{5}{2}\right)n^2+\left(4m-\frac{9}{2}\right)n+C^2\\
        &+(2mn-3n-2m^2-2m+6)C+(2m-3)m_2(G),
    \end{aligned}
    \end{equation*}
    where $C=|\{i\in[1,n]:2ms_i=0\}|$.
    \item If $m_{2m}(G)=0$, then
    \begin{equation*}
        \begin{aligned}
            M_2\leq (2m^2-4m+2)n^2+(4m-4)n+(2m-3)m_2(G).
        \end{aligned}
    \end{equation*}
\end{enumerate}
\end{lemma}
\begin{proof}
 Let $\bar{i},\bar{j}\in [1,n]\cup \{\infty\}$ such that $\bar{i}\equiv i \equiv l \pmod n$ and $\bar{j}\equiv j \equiv k \pmod n$. Since $(s_i,s_j)\in \Delta$, we have $i\not\equiv j \pmod n$, and so, $\bar{i}\neq \bar{j}$. Consider the following cases.
\begin{enumerate}
 \item   If $\bar{i}=\infty$ or $\bar{j}=\infty$, by exchanging the two sides of the equations, we assume that $\bar{i}=\infty$ and $\bar{j}\neq \infty$. Then $s_j$ and $s_k$ can be expressed as  $s_j=bs_{\bar{j}}$ and $s_k=cs_{\bar{j}}$, where $b,c \in [-m+1,m]^*$ such that $(b+c)s_{\bar{j}}=0$. 
    \begin{enumerate}
    \item If $b\neq m$ or $c\neq m$, without loss of generality, assume that $b\in [-m+1,m-1]^*$. Since $c s_{\bar{j}}=-b s_{\bar{j}}$ and $G\ge [-m+1,m]^* \diamond_2 S$, then $c=-b\in [-m+1,m-1]^*$. For each $\bar{j}$ and $b$, there is at most one equation satisfying the conditions. Then there are at most $2(m-1)n$ such equations. 
    \item If $b=c=m$, then $2m s_{\bar{j}}=0$. The number of such equations is at most $\left|\{i\in [1,n]: 2m s_i=0 \}\right|$.
    \end{enumerate}
\item If $\bar{i} \neq \infty$ and ${\bar{j}}\neq \infty$, let  $a,b,c,d \in [-m+1,m]^*$ such that  $s_i=a s_{\bar{i}}$, $s_j=b s_{\bar{j}}$,  $s_k=c s_{\bar{j}}$, and $s_l=d s_{\bar{i}}$. Then $(a+d)s_{\bar{i}}=(b+c)s_{\bar{j}}$. We proceed in the following subcases. 
    \begin{enumerate}
        \item If $a+d=0$ and $b+c=0$, for any $a,b\in [-m+1,m-1]^*$ and any distinct $\bar{i},\bar{j}\in [1,n]$, there is at most one equation. Therefore, there are at most $(2m-2)^2\binom{n}{2}=2(m-1)^2(n^2-n)$ such equations.

        \item If only one of $a+d$ and $b+c$ is $0$, we can assume without loss of generality that  $a+d=0$ and $b+c\neq 0$  by exchanging the two sides of the equations. Using a similar argument as in Case (1a), we can show that $b=c=m$, leading to $2m s_{\bar{j}}=0$. For a fixed $s_{\bar{j}}$, there are $2(m-1)$ choices for $(a,d)$ such that $a+d=0$, and $n-1$ choices for ${\bar{i}}$ such that $\bar{i}\neq \bar{j}$. Thus, there are at most $2(m-1)(n-1)\left|\{i\in [1,n]: 2m s_i=0 \}\right|$ such equations.

        \item If $a+d\neq 0$ and $b+c\neq 0$, then at least one of $a+d$ and $b+c$ must  equal  $2m$; otherwise, we can assume  without loss of generality that $d\neq m$ and $c\neq m$. Since $a s_{\bar{i}}-c s_{\bar{j}}=b s_{\bar{j}} -d s_{\bar{i}}$, this contradicts the condition $G\ge [-m+1,m]^* \diamond_2 S$. Now, assume that $a+d=2m$. Then $2ms_{\bar{i}}=(b+c)s_{\bar{j}}$. Hence, in this case there are most 
        \begin{align*} 
         & \sum_{\substack{x=-2(m-1) \\ x\neq 0}}^{2m-1} \tau(x) \left|\{(i,j)\in [1,n]^2:2m s_i=xs_j, i\neq j\}\right|\\
         & \ +  \left|\{i,j\}\subset [1,n]:2ms_i=2ms_j, i\neq j\}\right|
        \end{align*}
        such equations.
    \end{enumerate}
\end{enumerate}
    Therefore, we have that 
    \begin{equation}\label{Count:Equiv}
    \begin{aligned} 
        M_2\leq & \ 2(m-1)n+2(m-1)^2 (n^2-n)\\
        &+(2mn-2m-2n+3)\left|\{i\in [1,n]: 2m s_i=0 \}\right| \\
        &+\sum_{\substack{x=-2(m-1) \\ x\neq 0}}^{2m-1} \tau(x) \left|\{(i,j)\in [1,n]^2:2m s_i=xs_j, i\neq j\}\right|\\
        &+\left|\{i,j\}\subset [1,n]:2ms_i=2ms_j, i\neq j\}\right|
    \end{aligned}    
    \end{equation}
    such equations.
    
    In the following, we will estimate the right hand side of \eqref{Count:Equiv}.
    Notice that for each $i\in[1,n]$ satisfying $2m s_i=0$, we have $2ms_i\neq x s_j$ for any $x\in [-2m+2,2m-1]^*$ and $j\in [1,n]$, $j\neq i$. 
    
    For every odd $x \in [-2(m-1),2m-1]$, it can be written as $x=2q+1$, where $q\in [-m+1,m-1]$. Since $G\ge [-m+1,m]^* \diamond_2 S$, necessarily $(2q+1)s_j\neq 0$ for each $j \in [1,n]$. Then we claim that the set
    \begin{align*}
     & \ \{(i,j)\in[1,n]^2~:~  2m s_i=(2q+1) s_j,i\neq j\}
    \end{align*}
    contains at most  $n-C$ ordered pairs $(i,j)$. Otherwise, by the pigeonhole principle there exist two different pairs $(i,j)$ and $(i,j')$ such that $(2q+1)s_j=2m s_i=(2q+1)s_{j'}$, a contradiction. Then $\left|\{(i,j)\in[1,n]^2:2m s_i=(2q+1) s_j,i\neq j\}\right|\leq n-C$. We have
    \begin{align*}
        &\sum_{\scriptstyle x=-2(m-1) \atop\scriptstyle x\text{ odd}}^{2m-1} \tau(x) \left|\{(i,j)\in [1,n]^2:2m s_i=xs_j, i\neq j\}\right| \\
        &\leq (n-C)\sum_{\scriptstyle x=-2m+3 \atop\scriptstyle x\text{ odd}}^{2m-1} \tau(x)=(n-C)(2m^2-2m),
    \end{align*}
    where the equality holds due to Lemma~\ref{lemma:tau}.

    For every even $x\in [-2(m-1),2(m-1)]^*$, it can be written as $x=2q$, where $q\in [-m+1,m-1]^*$. If $2m s_i=2q s_j$ for $i\neq j\in[1,n]$, we have $2(ms_i-qs_j)=0$. Since $ms_i-qs_j\neq 0$, then $ord(ms_i-qs_j)=2$. Meanwhile, due to the splitting assumption, the elements from the set 
    \[\{ms_i- qs_j:i,j\in [1,n],q\in[-m+1,m-1]^*, i\neq j\}\cup \{ms_i:i\in[1,n]\}\] 
    are pairwise different. Therefore, for each $g\in G$ such that $ord(g)=2$, there exists at most one triple $(i,j,q)$ (or one $i$) such that $ms_i-q s_j=g$ (or $ms_i=g$). It follows that
    \[\sum_{\scriptstyle x=-2m+2 \atop\scriptstyle \text{even }x\neq 0 }^{2m-2}\left|\{(i,j)\in [1,n]^2:2m s_i=xs_j, i\neq j\}\right| +\left|\{i\in [1,n]: 2m s_i=0 \}\right| \leq m_2(G).\]
    Then
    \begin{align*}
        &\sum_{\scriptstyle x=-2m+2 \atop\scriptstyle x\text{ non-zero even}}^{2m-2}\tau(x) \left|\{(i,j)\in [1,n]^2:2m s_i=xs_j, i\neq j\}\right| \\
        &\leq (m_2(G)-C) \max_{\scriptstyle -2(m-1)\leq x \leq 2m-2 \atop\scriptstyle x\text{ non-zero even}}\tau(x)=(2m-3)(m_2(G)-C),
    \end{align*}
    where the equality holds since $\tau(x)$ attains its maximum value at $x=2$ when $x$ is a non-zero even number.
    
    Finally, we estimate $\left|\{i,j\}\subset [1,n]:2ms_i=2ms_j, i\neq j\}\right|$. As in the discussion above, we get $ord(s_i-s_j)=2m$. If $m_{2m}(G)=0$, this term is zero. If $m_{2m}(G)>0$, we have
    \begin{align*}
        &\left|\{i,j\}\subset [1,n]:2ms_i=2ms_j, i\neq j\}\right| \\
        &=\left|\{i,j\}\subset [1,n]:2ms_i=2ms_j=0, i\neq j\}\right| \\
        &+\left|\{i,j\}\subset [1,n]:2ms_i=2ms_j\neq 0, i\neq j\}\right|\\
        &\leq \binom{C}{2}+\binom{n-C}{2}=C^2-nC+\frac{n^2-n}{2},
    \end{align*}
    where we denote $\binom{1}{2}=\binom{0}{2}=0$.

    Therefore, continuing \eqref{Count:Equiv}, if $m_{2m}(G)>0$, we have that
    \begin{equation*}
    \begin{aligned} 
        M_2
        \leq & \ 2(m-1)n+2(m-1)^2 (n^2-n)+(2mn-2m-2n+3)C\\
        &+(n-C)(2m^2-2m)+(2m-3)(m_2(G)-C)+C^2-nC+\frac{n^2-n}{2}\\
        = &  \left(2m^2-4m+\frac{5}{2}\right)n^2+\left(4m-\frac{9}{2}\right)n+C^2\\
        &+(2mn-3n-2m^2-2m+6)C+(2m-3)m_2(G).
    \end{aligned}    
    \end{equation*}
    If $m_{2m}(G)=0$, then $C=0$, and so,
    \begin{equation*}
    \begin{aligned} 
        M_2 
        &\leq 2(m-1)n+2(m-1)^2 (n^2-n) +n(2m^2-2m)+(2m-3)m_2(G)\\
        &=(2m^2-4m+2)n^2+(4m-4)n+(2m-3)m_2(G).
    \end{aligned}    
    \end{equation*}
\end{proof}
\subsection{Proofs of the two lower bounds}\label{subsec:lbndpf}
\begin{proof}[Proof of Lemma~\ref{lm:lowerbound1}]
By Lemma~\ref{lm:M2}, we have 
\begin{equation*}
    \begin{aligned}
       M_2 \leq & \left(2m^2-4m+\frac{5}{2}\right)n^2+\left(4m-\frac{9}{2}\right)n+C^2\\
        &+(2mn-3n-2m^2-2m+6)C+(2m-3)m_2(G),
    \end{aligned}
    \end{equation*}
    when $m_{2m}(G)>0$; and 
\begin{equation*}
        \begin{aligned}
            M_2\leq (2m^2-4m+2)n^2+(4m-4)n+(2m-3)m_2(G),
        \end{aligned}
    \end{equation*}
    when $m_{2m}(G)=0$. Since $m\geq 4$ and $n\geq 4m$, we have $2mn-3n-2m^2-2m+6\geq 0$. Note that $C\leq m_2(G)$. Hence, in both cases, we have
    \begin{equation*}
    \begin{aligned}
       M_2 \leq & \left(2m^2-4m+\frac{5}{2}\right)n^2+\left(4m-\frac{9}{2}\right)n+m_2(G)^2\\
        &+(2mn-3n-2m^2+3)m_2(G).
    \end{aligned}
    \end{equation*}
    Combining this with \eqref{Inequality:Total} and \eqref{Maximum:Notequiv}, we obtain
    \begin{equation}\label{ineq:4}
        \begin{aligned}
            &(2m-1)^2\binom{n}{2}+(2m-1)n\\
            \leq & \ 2m^2 n+m_2(G)^2+(2mn-3n-2m^2+3)m_2(G)\\
             & +\left(2m^2-4m+\frac{5}{2}\right)n^2+\left(4m-\frac{9}{2}\right)n.
        \end{aligned}
    \end{equation}
    Let 
    $$f(x)\triangleq x^2-((2m-1)n-(2m^2-3))x+2m^2 n.$$
    Inequality \eqref{ineq:4} is equivalent to $f(n-m_2(G))\geq 0$. Note that $f(x)$ is decreasing when $x\leq n$. Since $n\geq 4m>2m+4$, we have $f(2m+\frac{3}{2})=4m^3+7m^2-2m^2n-\left(m-\frac{3}{2}\right)n-\frac{9}{4}<0$. Thus,
    $$n-m_2(G)\leq 2m+1,$$
    and so,
    $$m_2(G)+1\geq n-2m\geq \frac{n}{2}.$$
\end{proof}

In the proof of Lemma~\ref{lm:lowerbound1}, we used the simple estimate $C\leq m_2(G)$. To prove Lemma~\ref{lm:lowerbound2}, we require a more involved estimate of $C$.

\begin{lemma} \label{lemma:MinimumProjection} Let $G=G_1\times G_2$ be an Abelian group and $M=[-\km,\kp]^*$. Suppose that  $G=M\diamond_2 S$. Let   \[S_1\eqdef\{s_i\in S : \pi(s_i)\neq 0\},\]
where $\pi:G\to G_1$ is the natural projection.
Then
\[|S_1|> n\left(1-\frac{1}{\sqrt{|G_1|}}\right)-1.\]
\end{lemma}
\begin{proof}
    Denote $l\triangleq |S_1|$. If $G_1$ is trivial, the inequality obviously holds. Assume $G_1$ is nontrivial. For $g\in G_1\backslash \{0\}$, let
    $$f_1(g)\triangleq \left|\{(a,s_i):\pi(a s_i)=g,a\in M,s_i\in S_1\}\right|$$
    and
    $$f_2(g)\triangleq \left|\{(a,b,s_i,s_j):\pi(a s_i+b s_j)=g,a\text{ and }b\in M, s_i\neq s_j\in S_1, i<j\}\right|.$$
    Then
    $$\sum_{g\neq 0} f_1(g)\leq (\kp + \km)\binom{l}{1}$$
    and
    $$\sum_{g\neq 0} f_2(g)\leq (\kp + \km)^2\binom{l}{2}.$$
    
    Consider $|\pi^{-1}(g)|$. If $\pi(a s_i)=g$ for $a\in M$ and $s_i\in S$, we have $s_i\in S_1$. If $\pi(a s_i+b s_j)=g$ for $a,b\in M$ and $s_i,s_j\in S$, then at least one of $\set{s_i,s_j}$ belongs to  $S_1$. Therefore,     
    $$|G_2|=|\pi^{-1}(g)|=f_2(g)+f_1(g)+(\kp + \km)(n-l)f_1(g).$$
    It follows that
    \begin{align*}
     |G_2|(|G_1|-1) = & \ \sum_{g\in G_1 \backslash \set{0}} \abs{\pi^{-1}(g)} \\
    = & \ \sum_{g\in G_1 \backslash \set{0}} \parenv*{f_2(g)+f_1(g)+(\kp + \km)(n-l)f_1(g)} \\
    \leq & \ (\kp + \km)^2\binom{l}{2}+(1+(\kp + \km)(n-l))(\kp + \km)l.
    \end{align*}
    Rearranging the terms, we get
    \begin{equation}\label{Inequality:G_1}
        \begin{aligned}
            -\frac{(\kp + \km)^2}{2}l^2+\left(\left(n-\frac{1}{2}\right)(\kp + \km)^2+(\kp + \km)\right)l\ge |G|\left(1-\frac{1}{|G_1|}\right).
        \end{aligned}
    \end{equation}
    Denote $K\eqdef \kp+\km\geq 2$ and $q\eqdef |G_1|\geq 2$. Let 
    $$F(x)\eqdef -\frac{K^2}{2}x^2+\left(\left(n-\frac{1}{2}\right)K^2+K\right)x\ - \left(\frac{K^2}{2} n^2+\left(K-\frac{K^2}{2}\right)n+1\right) \left(1-\frac{1}{q}\right).$$
    Substituting $x=n\left(1-\frac{1}{\sqrt{q}}\right)-1$, we have 
    \begin{equation*}
        \begin{aligned}
            &F\left(n\left(1-\frac{1}{\sqrt{q}}\right)-1\right) \\
            =&-\frac{nK^2}{2}\left(\frac{1}{\sqrt{q}}+\frac{1}{q}\right)-nK\left(\frac{1}{\sqrt{q}}-\frac{1}{q}\right)- \frac{K^2}{2} -K  -\left(1-\frac{1}{q}\right) \\
            <& \ 0.
        \end{aligned}
    \end{equation*}
    Since $F(x)$ is  an increasing function when $x \leq n\left(1-\frac{1}{\sqrt{q}}\right)-1$, it follows from \eqref{Inequality:G_1} that 
    \begin{equation*} 
        \begin{aligned}
            l> n\left(1-\frac{1}{\sqrt{|G_1|}}\right)-1.
        \end{aligned}
    \end{equation*}
\end{proof}

\begin{lemma} \label{cor:maxC}
Let $G$ be an Abelian group and $S\subset G$ be a splitter set such that $G=[-m+1,m]^*\diamond_2 S$. Suppose that $G=G_1\times G_2$, where $gcd(|G_1|,2m)=1$.  Then   
    \begin{equation} \label{Maximum:C}
        C < \frac{n}{\sqrt{|G_1|}}+1.
    \end{equation}
\end{lemma}
\begin{proof} Let $\pi:G\to G_1$ be the natural projection.
    If $2m s_i=0$, then $2m\pi(s_i)=\pi(2ms_i)=0$. Since $\gcd(|G_1|,2m)=1$, necessarily $\pi(s_i)=0$. Hence,  we have  
    $$C\leq n-|S_1|  < \frac{n}{\sqrt{|G_1|}}+1,$$
    where the second inequality results from Lemma~\ref{lemma:MinimumProjection}.
\end{proof}

\begin{proof}[Proof of Lemma~\ref{lm:lowerbound2}]
    If $m_{2m}(G)>0$, it follows from \eqref{Inequality:Total}, together with Lemma~\ref{lm:M1} and Lemma~\ref{lm:M2}, that
    \begin{equation*}
        \begin{aligned}
        &2m^2 n+\left(2m^2-4m+\frac{5}{2}\right)n^2+\left(4m-\frac{9}{2}\right)n+C^2\\
        &+(2mn-3n-2m^2-2m+6)C+(2m-3)m_2(G) \\
        &\ge \frac{(2m-1)^2}{2}n^2-\frac{4m^2-8m+3}{2}n.
        \end{aligned}
    \end{equation*}
    Rearranging the terms, we get
    \begin{equation}\label{ineq:1}
        \begin{aligned}
        &C^2+(2mn-3n-2m^2-2m+6)C+(2m-3)m_2(G) \\
        &\ge (2m-2)n^2-(4m^2-3)n.
        \end{aligned}
    \end{equation}
    Let $q'\triangleq \max\{d\mid |G|:\gcd(2m,d)=1\}$. By Lemma~\ref{cor:maxC}, we have
    \begin{equation}\label{ineq:2}
    C < \frac{n}{\sqrt{q'}}+1.
    \end{equation}
    Note that when $n\ge m+4$ and $m\ge 2$, we have \begin{align*}
      &  \ 2mn-3n-2m^2-2m+6=(n-m)(2m-3)-5m+6\\
    \ge & \ 4(2m-3)-5m+6=3m-6 \geq 0.
    \end{align*}
    Combining the inequalities \eqref{ineq:1} and \eqref{ineq:2} we get
    \begin{equation*}
        \begin{aligned}
        &\left(\frac{n}{\sqrt{q'}}+1\right)^2+(2mn-3n-2m^2-2m+6)\left(\frac{n}{\sqrt{q'}}+1 \right)\\
        &+(2m-3)m_2(G)> (2m-2)n^2-(4m^2-3)n.
        \end{aligned}
    \end{equation*}
    If $m_{2m}(G)=0$, combining Lemma~\ref{lm:M1} and Lemma~\ref{lm:M2} we have
    $$(2m-3)m_2(G)\geq \left(2m-\frac{5}{2}\right)n^2-\left(4m^2-\frac{5}{2}\right)n.$$
\end{proof}

\section{Upper bounds on $m_2(G)$}\label{Sec:Upper}
In this section, we are going to prove the following two upper bounds on $m_2(G)$.

\begin{lemma}\label{lm:upperbound1}
Suppose that $G = [-m+1,m]^* \splt_{2} S$. Let $s\triangleq\lfloor \log_2 m\rfloor$. Then,
$$\log_2 (m_2(G)+1) \leq \frac{2}{s+1} \log_2 n+\lfloor\log_2 s\rfloor+\frac{s-2^{\lfloor\log_2 s\rfloor+1}}{s+1}+5.$$
\end{lemma}

\begin{lemma}\label{lm:upperbound2}
    Suppose that $G\ge [-m+1,m]^*\diamond_2 S$ and $|G|=2^a\times q$, where $q$ is an odd integer. Then
    \begin{equation} \label{Maximum:m2}
        m_2(G)\leq \frac{|G|}{\sqrt{q\left( (m-1)^2\frac{n^2-n}{2}+(m-1)n+1\right)}}-1.
    \end{equation}
\end{lemma}

The second lemma is sharper when $m=2,3$, and accounts for the effect of the odd factor $q$ of $|G|$.

The proofs of both lemmas make use of the following auxiliary result.
\begin{lemma}\label{lemma:q_times}
    Let $G=\prod_l \Z_l^{d_l}$, where $l$ runs over all prime power. Let $q$ be a prime power and assume that $G=\prod_{l\mid q} \Z_l^{d_l} \times G'$. If $G \geq [-\km, \kp]^* \splt_t S$ for some splitter set $S$, then $G'\geq [-\lfloor \frac{\km}{q} \rfloor, \lfloor \frac{\kp}{q} \rfloor]^* \splt_t S'$ for some splitter set $S'$ with $|S'|=|S|$.
\end{lemma}
\begin{proof}
    Since $G=\prod_{l|q} \Z_l^{d_l} \times G'$, for each $s_i \in S$ and $l\mid q$, the image of the natural projection of $q s_i$ to $\Z_l$ is $0$. Let  $g_i$ denote the image of the natural projection of $qs_i$ to $G'$. Denote $S'\triangleq \{g_1,g_2,\ldots,g_n\}\subseteq G'$. Since $G\ge [-\km,\kp]^* \splt_t S$, it is straightforward to verify that  $G'\geq [-\lfloor \frac{\km}{q} \rfloor, \lfloor \frac{\kp}{q} \rfloor]^*\splt_t S'$.
\end{proof}

\subsection{The first upper bound}
We get a series of constraints on $d_{2^i}$ by applying Lemma~\ref{lemma:q_times} with even $q$. 
\begin{lemma}\label{lemma:constraints_d2i}
    Let $G=\prod_l \Z_l^{d_l}$, where $l$ runs over all prime power. Suppose that $G =[-m+1, m]^* \splt_2 S$. Then,
    $$\sum_{i=1}^j i d_{2^i} \leq 2j+2.$$
    for all $1 \leq j \leq \lfloor \log_2 m\rfloor$.
\end{lemma}

\begin{proof}
    Applying Lemma~\ref{lemma:q_times} with $q = 2^j$, we have that
    $$G' \geq [-\lfloor \frac{m-1}{2^j} \rfloor, \lfloor \frac{m}{2^j} \rfloor] \splt_2 S'.$$
    Then,
    $$\prod_{i=1}^j 2^{i d_{2^i}} = \frac{|G|}{|G'|} \leq \frac{f(2m-1)}{f(\lfloor \frac{m-1}{q} \rfloor + \lfloor \frac{m}{q} \rfloor)} \leq q^{2} \left( \frac{\frac{m-1}{q} + \frac{m}{q}}{\lfloor \frac{m-1}{q} \rfloor + \lfloor \frac{m}{q} \rfloor} \right)^{2},$$
    where
    $$f(m)=\binom{n}{2} x^2+n x+1$$ and the last inequality holds as $2m-1 \geq \lfloor \frac{m-1}{q}\rfloor+\lfloor\frac{m}{q}\rfloor$.
    
    Let $x \triangleq \frac{m-1}{q}$ and $y \triangleq \frac{m}{q}$. If $\lfloor x\rfloor=0$ and $\lfloor y \rfloor=1$, we have $m=q$, then $y=1$, 
    $$\frac{x+y}{\floor{x}+\floor{y}}\leq 2y= 2.$$
    If $\floor{x}\geq 1$,
    $$\frac{x+y}{\floor{x}+\floor{y}}=1+\frac{x-\lfloor x \rfloor + y - \lfloor y \rfloor}{\lfloor x \rfloor + \lfloor y \rfloor}\leq 2.$$
    Therefore,
    $$\prod_{i=1}^j 2^{i d_{2^i}}\leq 4q^{2}=2^{2j+2}.$$
    Taking $\log_{2}$ on both sides, we obtain
    $$\sum_{i=1}^j i d_{2^i} \leq 2 j+2.$$
\end{proof} 

\begin{lemma}\label{lm:techlem}
     Let $x_1,\ldots,x_s$ be non-negative integers with $s\geq 2$. Suppose that for each $1\leq j\leq s$, 
    $$\sum_{i=1}^j i x_i\leq 2j+2.$$
    Then we have
    $$\sum_{i=1}^s (s+1-i)x_i\leq 4s+r(s+1)-2^{r+1}+2,$$
    where $r=\lfloor \log_2 s \rfloor$.
\end{lemma}
\begin{proof} Denote 
    $$f(x_1,x_2,\ldots,x_s)=\sum_{i=1}^s (s+1-i)x_i.$$
    In the following we will show that under the constraints $\sum_{i=1}^j i x_i\leq 2j+2$ for all $1\leq j\leq s$, the objective function $f(x_1,x_2,\ldots,x_s)$  attains its maximum if and only if the variables take the following values:
    \begin{equation}\label{eq:solution}
    x_i=
    \begin{cases}
        4, &\text{if }i=1,\\
        1, &\text{if }i=2^l,\ 1\leq l\leq r,\\
        0, &\text{otherwise}.
    \end{cases}
    \end{equation}
    The conclusion then follows, as in this case, we obtain
    $$f(x_1,x_2,\ldots,x_s)=4s+\sum_{i=1}^r (s+1-2^i)=4s+r(s+1)-2^{r+1}-2.$$

    First, assume $x_1<4$.  Let $k\geq 2$ be the smallest integer such that $x_k>0$ and set
    $$x_1'=x_1+1, \quad x_k'=x_k-1, \quad x'_i=x_i \textrm{ for all } i\neq 1,k.$$ If no such $k$ exists, then $x_i=0$ for all $i\geq 2$, and we instead set 
    $$x_1'=x_1+1,\quad x'_i=x_i \textrm{ for all } i\geq 2.$$
    It is easy to verify that all constraints still hold and $$f(x'_1,x'_2,\ldots, x'_s)-f(x_1,x_2,\ldots,x_s)=k-1>0.$$
    Now, we consider the case of $x_1=4$. For any $2\leq i \leq s$, the constraints imply that
    $$2x_2+\cdots+i x_i\leq 2i-2.$$
    Thus we must have $x_i\leq 1$. 
    Suppose that there are non-power-of-2 indices $i$'s where $x_i=1$. Let $k$ be the smallest such index and denote $l=\lfloor \log_2 k\rfloor$. Then we have $2^l<k<2^{l+1}$. If $x_{2^i}=1$ for all $1\leq i\leq l$, then
    $$\sum_{i=1}^k j x_j=4+2x_2+\cdots+2^l x_{2^l}+k x_k=2^{l+1}+2+k>2k+2,$$
    which contradicts the constraint. Thus, we choose an arbitrary $i_0\in [1,l]$ such that $x_{2^{i_0}}=0$. Let
    $$x'_{2^{i_0}}=1,\quad x'_k=0, \quad x'_j=x_j \textrm{ for all } j\neq \set{1,i_0,k}.$$
    It is easy to verify that all constraints still hold, and $$f(x'_1,x'_2,\ldots, x'_s)-f(x_1,x_2,\ldots,x_s)=k-2^i>0.$$
    The argument above shows that the object function attains its maximum only when $x_1=4$ and $x_i=0$ for all indices $i$  that are not powers of 2. A simple calculation further shows that the unique optimal solution is given by \eqref{eq:solution}.
\end{proof}

\begin{proof}[Proof of Lemma~\ref{lm:upperbound1}]
 Let $G=\prod_l \Z_l^{d_l}$, where $l$ runs over all prime power. Then $m_2(G)+1=2^{\sum_{i=1}^\infty d_{2^i}}$. Hence,
    \begin{equation*}
        \begin{aligned}
            &(s+1) \log_2 (m_2(G)+1)-\log_2 |G|\\
            \leq &(s+1)\sum_{i=1}^{\infty}d_{2^i}-\sum_{i=1}^{\infty}i d_{2^i} \\
            \leq &\sum_{i=1}^{s} (s+1-i) d_{2^i} \\
            \leq & 4s+r(s+1)-2^{r+1}+2\\
            =&4s+(s+1)\lfloor\log_2 s\rfloor-2^{\lfloor\log_2 s\rfloor+1}+2,
        \end{aligned}
    \end{equation*}
    where the last inequality  follows from Lemma~\ref{lemma:constraints_d2i} and Lemma~\ref{lm:techlem}.
    
    Since $|G| = (2m-1)^2 \binom{n}{2}+(2m-1)n+1 \leq 2m^2 n^2$, taking the logarithm base $2$ on both sides gives
    $$\log_2 |G| \leq 2 \log_2 n+ 2\log_2 m +1.$$    
    Thus,
    \begin{equation*}
        \begin{aligned}
            \log_2 (m_2(G) + 1) 
            &\leq \frac{1}{s+1} \left(\log_2 |G|+4s+(s+1)\lfloor\log_2 s\rfloor-2^{\lfloor\log_2 s\rfloor+1}+2 \right)\\
            &\leq \lfloor\log_2 s\rfloor+\frac{1}{s+1} \left(2\log_2 n+2\log_2 m +4s-2^{\lfloor\log_2 s\rfloor+1}+3\right)\\
            &<\lfloor\log_2 s\rfloor+\frac{1}{s+1}\left(2\log_2 n+2s+2+4s-2^{\lfloor\log_2 s\rfloor+1}+3\right)\\
            &=\frac{2}{s+1} \log_2 n+\lfloor\log_2 s\rfloor+\frac{s-2^{\lfloor\log_2 s\rfloor+1}}{s+1}+5.
        \end{aligned}
    \end{equation*}

\end{proof}

\subsection{The second upper bound}

In the proof of Lemma~\ref{lm:upperbound1}, the odd factors of $G$ were ignored. In this subsection, we take these factors into account to derive more precise bounds on $m_2(G)$, particularly for  the cases $m=2,3$.

\begin{proof}[Proof of Lemma~\ref{lm:upperbound2}]
    Suppose that $G=\prod_{i\geq 1}\Z_{2^i}^{d_{2^i}}\times G''$, where $|G''|=q$. Denote $k\triangleq\sum_{i=2}^{\infty}d_{2^i}$
    Then 
    \[|G|=\prod_{i\geq 1}2^{id_{2^i}}\times q \ge 2^{d_2+2k} \times q.\]
    It follows that $2^{d_2+k} \leq \sqrt{2^{d_2} \abs{G}/q}.$
    Hence, 
    \begin{equation}\label{eq:m2bnd}
    m_2(G)=2^{d_2+k}-1\leq \sqrt{\frac{2^{d_2} |G|}{q}}-1.
    \end{equation}
    Applying Lemma~\ref{lemma:q_times} with ``$q=2$", we have $G'\ge [-\lfloor\frac{m-1}{2}\rfloor, \lfloor\frac{m}{2}\rfloor]^*\diamond_2 S'$. Therefore
    $$\frac{|G|}{2^{d_2}}=|G'|\ge (m-1)^2\frac{n^2-n}{2}+(m-1)n+1.$$
    Rearranging the terms, we get
   \begin{equation}\label{eq:2alpha}
   2^{d_2} \leq \frac{\abs{G}}{ (m-1)^2\frac{n^2-n}{2}+(m-1)n+1}. 
   \end{equation}
   Substituting \eqref{eq:2alpha} in \eqref{eq:m2bnd} gives 
    \begin{equation*} 
        m_2(G)\leq \frac{|G|}{\sqrt{q\left( (m-1)^2\frac{n^2-n}{2}+(m-1)n+1\right)}}-1.
    \end{equation*}
\end{proof}

\section{Proof of the main result}\label{Sec:main}
In this section, we are going to prove 
Theorem~\ref{thm:main}, which is established by proving the following two results, supplemented by a computer search.

\begin{proposition}\label{theorem1}
    Let $m\geq 4$ and $n\geq 4m$. If $\cB(n,2,m,m-1)$ lattice-tiles $\mathbb{Z}^n$, then one of the following cases holds:
    \begin{enumerate}
        \item $4\leq m \leq 7$ and $n\leq 524288$;
        \item $8\leq m \leq 15$ and $n\leq 11585$;
        \item $16\leq m \leq 31$ and $n\leq 4096$;
        \item $32\leq m \leq 63$ and $n\leq 2435$;
        \item $64\leq m \leq 127$ and $n\leq 1782$;
        \item $128\leq m \leq 255$ and $n\leq 1448$;
        \item $256\leq m \leq 511$ and $n\leq 1378$.
    \end{enumerate}
\end{proposition}

\begin{proof}
    Combining Lemma~\ref{lm:lowerbound1} and Lemma~\ref{lm:upperbound1}, we have
    $$\log_2 n -1\leq\log_2(m_2(G)+1)\leq\frac{2}{s+1} \log_2 n+\lfloor\log_2 s\rfloor+\frac{s-2^{\lfloor\log_2 s\rfloor+1}}{s+1}+5.$$
    Therefore,
    $$s+2\leq\log_2 n\leq \frac{s+1}{s-1}(\lfloor\log_2 s\rfloor+\frac{s-2^{\lfloor\log_2 s\rfloor+1}}{s+1}+6),$$
    where the first inequality holds as $n\geq 4m\geq 2^{s+2}$.
    It follows that $s\leq 8$, and for each $s$, an upper bound for $n$,  as stated in Proposition~\ref{theorem1}, can be determined. 
\end{proof}


\begin{proposition}\label{theorem2}
    Let $n\ge 3$ and $m\ge 2$. If $\cB(n,2,m,m-1)$ lattice-tiles $\mathbb{Z}^n$, then  one of the following conditions must hold:
    \begin{enumerate}
        \item $n\leq (2+2\sqrt{2})m+2 $;
        \item $4m \mid n^2-3n+2$, and for every prime number $p\mid (2m-1)^2\frac{n^2-n}{2}+(2m-1)n+1$, $p$ must divide $2m$;
        \item $4m \mid n^2-3n+2$ and $n \leq \left(\frac{5+\sqrt{3}}{2}+\sqrt{2}+\sqrt{6}\right)m+4$.
    \end{enumerate}
\end{proposition}

\begin{proof}
    Since $\cB(n,2,m,m-1)$ lattice-tiles $\mathbb{Z}^n$, there is an Abelian group $G$ with $\abs{G}=\abs{\cB(n,2,m,m-1)}$ such that $G=[-m+1,m]^* \diamond_2 S$. We first examine the case where $m_{2m}(G)=0$. 
    By Lemma~\ref{lm:lowerbound2} and Lemma~\ref{lm:upperbound2}, we have
    \begin{equation*}
        \begin{aligned}
            &(2m-3)\left(\frac{(2m-1)^2 \frac{n^2-n}{2}+(2m-1)n+1}{\sqrt{q\left( (m-1)^2\frac{n^2-n}{2}+(m-1)n+1\right)}}-1\right) \\
            &\ge \left(2m-\frac{3}{2}\right)n^2-\left(4m^2-\frac{5}{2}\right)n.
        \end{aligned}
    \end{equation*}
   We are going to show that $n\leq (2+\sqrt{2})m+2$. Since 
    $$\left(2m-\frac{3}{2}\right)n^2-\left(4m^2-\frac{5}{2}\right)n\ge n\left(2m-\frac{3}{2}\right)(n-2m-2),$$
    $$\sqrt{\left( (m-1)^2\frac{n^2-n}{2}+(m-1)n+1\right)}\ge \frac{(m-1)(n-1)}{\sqrt{2}},$$
    and
    $$(2m-1)^2 \frac{n^2-n}{2}+(2m-1)n+1\leq \frac{(2m-1)^2}{2}n^2,$$
    we have
    $$(2m-3)\frac{\sqrt{2}(2m-1)^2n^2}{2\sqrt{q}(m-1)(n-1)} \ge n\left(2m-\frac{3}{2}\right)(n-2m-2).$$
    Suppose that $n\ge 2m+2$, then $\frac{n}{n-1}\leq \frac{2m+2}{2m+1}$, it follows that
    \[ (2m-3)\frac{\sqrt{2}(2m-1)^2(2m+2) }{2\sqrt{q}(m-1)(2m+1)} \ge \left(2m-\frac{3}{2}\right)(n-2m-2).\]
    Rearranging the terms,
    we get
    \begin{equation} \label{Inequality:nmid}
        \begin{aligned}
            &n\leq 2m+2+\frac{\sqrt{2}(2m-1)^2(2m-3)(2m+2)}{2\sqrt{q}(2m-3/2)(m-1)(2m+1)}\\
            &= 2m+2+\frac{\sqrt{2}}{\sqrt{q}}\cdot \frac{(2m-1)^2}{2m-3/2} \cdot \frac{(2m-3)(2m+2)}{(2m-2)(2m+1)}\\
            &= 2m+2+\frac{\sqrt{2}}{\sqrt{q}}\cdot \left(2m-\frac{m-1}{2m-3/2}\right) \cdot \frac{4m^2-2m-6}{4m^2-2m-2}\\
            &< \parenv*{2+\frac{2\sqrt{2}}{\sqrt{q}}}m+2 \\
            &\leq (2+2\sqrt{2})m+2.
        \end{aligned}
    \end{equation}
    
    Now, we proceed with the case where $m_{2m}(G)>0$. Here,  $2m\mid \abs{G}$. Given that $|G|=(2m-1)^2\frac{n^2-n}{2}+(2m-1)n+1\equiv \frac{n^2-3n+2}{2} \pmod{2m}$, it follows that $4m \mid n^2-3n+2$. 
    
    Suppose that $|G|=2^a \times q$, where $q$ is an odd integer. By Lemma~\ref{lm:upperbound2}, we have
    \begin{equation}\label{ieq-2}
    m_2(G)\leq \frac{|G|}{\sqrt{q\left( (m-1)^2\frac{n^2-n}{2}+(m-1)n+1\right)}}-1.
    \end{equation}
    Let $q'=\max\{d\mid |G|:\gcd(2m,d)=1\}$. By Lemma~\ref{lm:lowerbound2}, we have
    \begin{equation*}
    \begin{aligned}
        &\left(\frac{n}{\sqrt{q'}}+1\right)^2+(2mn-3n-2m^2-2m+6)\left(\frac{n}{\sqrt{q'}}+1 \right)\\
        &+(2m-3)m_2(G)>(2m-2)n^2-(4m^2-3)n.
    \end{aligned}
    \end{equation*}
    Combining the inequality \eqref{ieq-2}, we get
    \begin{equation*}
        \begin{aligned}
        &\left(\frac{n}{\sqrt{q'}}+1\right)^2+(2mn-3n-2m^2-2m+6)\left(\frac{n}{\sqrt{q'}}+1 \right)\\
        &+(2m-3)\left(\frac{(2m-1)^2 \frac{n^2-n}{2}+(2m-1)n+1}{\sqrt{q\left( (m-1)^2\frac{n^2-n}{2}+(m-1)n+1\right)}}-1\right) \\
        &> (2m-2)n^2-(4m^2-3)n.
        \end{aligned}
    \end{equation*}
    
    Rearranging it we get
    
    \begin{equation}\label{Inequality:mid}
        \begin{aligned}
        &(2m-3)\left(\frac{(2m-1)^2 \frac{n^2-n}{2}+(2m-1)n+1}{\sqrt{q\left((m-1)^2\frac{n^2-n}{2}+(m-1)n+1\right)}}-1\right)-2m^2-2m+7 \\
        &> \left(2m-2-\frac{1}{q'}-\frac{2m-3}{\sqrt{q'}}\right)n^2-\left(4m^2-3+\frac{2}{\sqrt{q'}}-\frac{2m^2+2m-6}{\sqrt{q'}}+2m-3\right)n.
        \end{aligned}
    \end{equation}
    We claim that
    $$LHS\text{ of (\ref{Inequality:mid})}<\frac{\sqrt{2}}{\sqrt{q}}(2m-1)^2 n.$$
    If $m=2$, 
    $$\frac{9\frac{n^2-n}{2}+3n+1}{\sqrt{q\left(\frac{n^2-n}{2}+n+1\right)}}-1-5<\frac{9n^2}{2\sqrt{q}\frac{n}{\sqrt{2}}}<\frac{9\sqrt{2}}{\sqrt{q}}n.$$
    If $m\geq 3$, noting that $n\geq 3$, we have that
    \begin{equation*}
        \begin{aligned}
             &  (m-1)^2\frac{n^2-n}{2}+(m-1)n+1-\frac{(m-1)^2(n-4/7)^2}{2} \\
            =&  \left(\frac{(m-1)^2}{14}+m-1\right)n-(m-1)^2\frac{8}{49}+1\geq 0,
        \end{aligned}
    \end{equation*}
    and
    \begin{equation*}
        \begin{aligned}
            &(2m-1)^2 \frac{n^2-n}{2}+(2m-1)n+1 - \frac{(2m-1)^2}{2}\left(n^2-\frac{4}{7}n\right) \\
            = &\left(-\frac{(2m-1)^2}{2}\cdot \frac{3}{7}+2m-1\right)n+1 \leq \left(-\frac{25}{2}\cdot\frac{3}{7}+5\right)n+1 \\
            =&-\frac{5}{14}n+1\leq 0.
        \end{aligned}
    \end{equation*}
    It follows that
    $$\sqrt{(m-1)^2\frac{n^2-n}{2}+(m-1)n+1}\ge \frac{(m-1)\left(n-\frac{4}{7}\right)}{\sqrt{2}}$$
    and
    $$(2m-1)^2 \frac{n^2-n}{2}+(2m-1)n+1\leq \frac{(2m-1)^2}{2}\left(n^2-\frac{4}{7}n\right).$$
    Hence,
    $$LHS\text{ of (\ref{Inequality:mid})}\leq \frac{\sqrt{2}(2m-3)(2m-1)^2\left(n^2-\frac{4}{7}n\right)}{2\sqrt{q}(m-1)\left(n-\frac{4}{7}\right)}<\frac{\sqrt{2}}{\sqrt{q}}(2m-1)^2 n.$$
    Therefore,
    \begin{equation}\label{Inequality:mid_2}
        \begin{aligned}
            \frac{\sqrt{2}}{\sqrt{q}}(2m-1)^2 n \geq &\left(2m-2-\frac{1}{q'}-\frac{2m-3}{\sqrt{q'}}\right)n^2 \\
            &-\left(4m^2-3-\frac{2m^2+2m-8}{\sqrt{q'}}+2m-3\right)n.
        \end{aligned}
    \end{equation}

    If $q'=1$, then all nontrivial factors of $|G|$ divides $2m$. In other words,  for every prime number $p\mid (2m-1)^2\frac{n^2-n}{2}+(2m-1)n+1$, we have that $p\mid 2m$. Now, suppose that $q'\ge 3$, then $q\ge q'\ge 3$. Hence, 
    $$LHS\text{ of (\ref{Inequality:mid_2})}< \frac{\sqrt{2}}{\sqrt{3}}(2m-1)^2 n.$$
    Since $RHS$ of (\ref{Inequality:mid}) is an increasing function with respect to $q'$,    
    
    \begin{equation*}
        \begin{aligned}
            RHS\text{ of (\ref{Inequality:mid_2})}\geq & \left(2m-2-\frac{1}{3}-\frac{2m-3}{\sqrt{3}}\right)n^2\\
            &-\left(4m^2-3+\frac{2}{\sqrt{3}}-\frac{2m^2+2m-6}{\sqrt{3}}+2m-3\right)n.
        \end{aligned}
    \end{equation*} 
    The inequality (\ref{Inequality:mid_2}) yields
    \begin{equation*}
        \begin{aligned}
            n&\leq \frac{4m^2-3+\frac{2}{\sqrt{3}}-\frac{2m^2+2m-6}{\sqrt{3}}+2m-3+\frac{\sqrt{2}}{\sqrt{3}}(2m-1)^2}{2m-\frac{7}{3}-\frac{2m-3}{\sqrt{3}}} \\
            &=\frac{\left(4+\frac{4\sqrt{2}-2}{\sqrt{3}}\right)m^2-\left(\frac{4\sqrt{2}+2}{\sqrt{3}}{-2}\right)m-{6}+2\sqrt{3}+\frac{\sqrt{2}{+2}}{\sqrt{3}}}{\left(2-\frac{2}{\sqrt{3}}\right)m-\frac{7}{3}+\sqrt{3}} \\
            &\leq \left(\frac{5+\sqrt{3}}{2}+\sqrt{2}+\sqrt{6}\right)m+{4}.
        \end{aligned}
    \end{equation*}
    
\end{proof}

Now, we can prove Theorem~\ref{thm:main}.


\begin{proof}[Proof of Theorem~\ref{thm:main}] We fist consider the case of $m=2$.
    Suppose that $\cB(n,2,2,1)$ lattice-tiles $\mathbb{Z}^n$, then Proposition~\ref{theorem2} shows that one of the following holds:
    \begin{enumerate}
        \item[(A1)] $n\leq 2(2+2\sqrt{2})+2$;
        \item[(A2)] $8 \mid n^2-3n+2$, and $\frac{9}{2}n^2-\frac{3}{2}n+1$ is a power of 2;
        \item[(A3)] $8 \mid n^2-3n+2$, $n \leq 2\left(\frac{5+\sqrt{3}}{2}+\sqrt{2}+\sqrt{6}\right)+4$.
    \end{enumerate}
    
    For (A1), we have $n\leq 11$. A computer search rules out the values $n \in [3, 11]$.

    For (A2), assume that $\frac{9}{2}n^2-\frac{3}{2}n+1=2^\alpha$ for some integer $\alpha$. Rearranging this we get $$(6n-1)^2=2^{\alpha+3}-7.$$ 
    By Lemma~\ref{lemma:Nagell}, we have $6n-1\in \set{1,3,5,11,181}$. So, $n\in \set{\frac{1}{3},\frac{2}{3},1,2,\frac{91}{3}}$, a contradiction. 

    For (A3), given that $8\mid (n-1)(n-2)$ and $n\le 18$, the possible pairs of $(n,|G|)$ are $(9,352=2^5\times 11),(10,436=2^2\times 109),(17,1276=2^2\times 319),(18,1432=2^3\times 179)$. Using \eqref{Inequality:mid_2}, we obtain a contradiction. 


    Next, we consider the case of $m=3$.
    Suppose that $\cB(n,2,3,2)$ lattice-tiles $\mathbb{Z}^n$. Then, similarly, by  Proposition~\ref{theorem2}, one of the following holds:
    \begin{enumerate}
        \item[(B1)] $n\leq 3(2+2\sqrt{2})+2$;
        \item[(B2)] $12 \mid n^2-3n+2$, and all prime factors of $\frac{25}{2}n^2-\frac{15}{2}n+1$ are 2 and 3;
        \item[(B3)] $12 \mid n^2-3n+2$, $n \leq 3\left(\frac{5+\sqrt{3}}{2}+\sqrt{2}+\sqrt{6}\right)+4$.
    \end{enumerate}
    
    For (B1), we have $n\leq 16$. A computer search rules out these values.  
    
    For (B2), assume that $\frac{25}{2}n^2-\frac{15}{2}n+1=2^\alpha 3^\beta$ for some integer $\alpha$ and $\beta$. Rearranging this we get $$(5n-1)(5n-2)=2^{\alpha+1}3^\beta.$$ 
    Then $(5n-1,5n-2)\in\set{(2^{\alpha+1},3^\beta), (3^\beta, 2^{\alpha+1})}$. 
    If $5n-1=2^{\alpha+1}$ and $5n-2=3^\beta$, we have $\beta \equiv 1 \pmod 4$. Let $\beta=4k+1$, then $2^{\alpha+1}=3^\beta+1= 3^{4k+1}+1\equiv3\times 9^{2k}+1\equiv 4 \pmod 8$. Therefore $2^{\alpha+1}=4$, and so, $n=1$, a contradiction. If $5n-1=3^\beta$ and $5n-2=2^{\alpha+1}$, then $\beta \equiv 2 \pmod 4$. Let $\beta=4k+2$, then we have $3^{4k+2}=2^{\alpha+1}+1$. Therefore $(3^{2k+1}+1)(3^{2k+1}-1)=2^{\alpha+1}$. Since $\gcd(3^{2k+1}-1,3^{2k+1}+1)=1$ or $2$, we have $3^{2k+1}=3$, and so, $n=2$, a contradiction. 

    For (B3), given that $12\mid (n-1)(n-2)$ and $n\le 27$, the possible pairs of $(n,|G|)$ are $(17,3486), (22,5886),(25,7626), (26,8256)$. Using \eqref{Inequality:mid_2}, we obtain a contradiction in these cases. 

    Finally, for $m\geq 4$, a quick computer search\footnote{Note that we only need to check that when $4m \mid n^2-3n+2$, the size of the ball has a prime factor other than  $2m$. This can be done quickly using the Euclidean algorithm.} confirms that none of the values of $m$ and $n$ listed in Proposition~\ref{theorem1} satisfy condition (2) of Proposition~\ref{theorem2}. The conclusion then follows by combining Proposition~\ref{theorem1} and Proposition~\ref{theorem2}.
\end{proof}


\begin{thebibliography}{19}

\bibitem{buzaglo2012tilings}
S.~Buzaglo and T.~Etzion.
\newblock Tilings with $ n $-dimensional chairs and their applications to
  asymmetric codes.
\newblock {\em IEEE Transactions on Information Theory}, 59(3):1573--1582,
  2012.

\bibitem{BuzEtz13}
S.~Buzaglo and T.~Etzion.
\newblock Tilings by $(0.5,n)$-crosses and perfect codes.
\newblock {\em SIAM Journal on Discrete Mathematics}, 27(2):1067--1081, 2013.

\bibitem{cassuto2010codes}
Y.~Cassuto, M.~Schwartz, V.~Bohossian, and J.~Bruck.
\newblock Codes for asymmetric limited-magnitude errors with application to
  multilevel flash memories.
\newblock {\em IEEE Transactions on Information Theory}, 56(4):1582--1595,
  2010.

\bibitem{hamaker1984combinatorial}
W.~Hamaker and S.~Stein.
\newblock Combinatorial packings of $\mathbb{R}^3$ by certain error spheres.
\newblock {\em IEEE Transactions on Information Theory}, 30(2):364--368, 1984.

\bibitem{hickerson1986abelian}
D.~Hickerson and S.~Stein.
\newblock Abelian groups and packing by semicrosses.
\newblock {\em Pacific Journal of Mathematics}, 122(1):95--109, 1986.

\bibitem{Jainetal2020ISIT}
S.~Jain, F.~Farnoud, M.~Schwartz, and J.~Bruck.
\newblock Coding for optimized writing rate in {DNA} storage.
\newblock In {\em Proceedings of the 2020 IEEE International Symposium on
  Information Theory (ISIT2020), Los Angeles, CA, USA}, pages 711--716, 2020.

\bibitem{klove2011some}
T.~Kl{\o}ve, J.~Luo, I.~Naydenova, and S.~Yari.
\newblock Some codes correcting asymmetric errors of limited magnitude.
\newblock {\em IEEE Transactions on Information Theory}, 57(11):7459--7472,
  2011.

\bibitem{LeeKalGoeBolChur19}
H.~H. Lee, R.~Kalhor, N.~Goela, J.~Bolot, and G.~M. Church.
\newblock Terminator-free template-independent enzymatic {DNA} synthesis for
  digital information storage.
\newblock {\em Nature Communications}, 10(2383):1--12, 2019.

\bibitem{nagell1961diophantine}
T.~Nagell.
\newblock The diophantine equation $x^2+7=2^n$.
\newblock {\em Arkiv f{\"o}r Matematik}, 4(2):185--187, 1961.

\bibitem{schwartz2011quasi}
M.~Schwartz.
\newblock Quasi-cross lattice tilings with applications to flash memory.
\newblock {\em IEEE Transactions on Information Theory}, 58(4):2397--2405,
  2011.

\bibitem{schwartz2014non}
M.~Schwartz.
\newblock On the non-existence of lattice tilings by quasi-crosses.
\newblock {\em European Journal of Combinatorics}, 36:130--142, 2014.

\bibitem{stein1984packings}
S.~Stein.
\newblock Packings of ${R}^n$ by certain error spheres.
\newblock {\em IEEE Transactions on Information Theory}, 30(2):356--363, 1984.

\bibitem{stein1990notched}
S.~Stein.
\newblock The notched cube tiles ${R}^n$.
\newblock {\em Discrete mathematics}, 80(3):335--337, 1990.

\bibitem{stein1994algebra}
S.~Stein and S.~Szab{\'o}.
\newblock {\em Algebra and Tiling: Homomorphisms in the service of Geometry}.
\newblock Number~25. Cambridge University Press, 1994.

\bibitem{WeiSch21}
H.~Wei and M.~Schwartz.
\newblock Improved coding over sets for {DNA}-based data storage.
\newblock {\em IEEE Transactions on Information Theory}, 68(1):118--129, 2021.

\bibitem{WEI2022103450}
H.~Wei and M.~Schwartz.
\newblock On tilings of asymmetric limited-magnitude balls.
\newblock {\em European Journal of Combinatorics}, 100:103450, 2022.

\bibitem{wei2021lattice}
H.~Wei, X.~Wang, and M.~Schwartz.
\newblock On lattice packings and coverings of asymmetric limited-magnitude
  balls.
\newblock {\em IEEE Transactions on Information Theory}, 67(8):5104--5115,
  2021.

\bibitem{yari2013some}
S.~Yari, T.~Kl{\o}ve, and B.~Bose.
\newblock Some codes correcting unbalanced errors of limited magnitude for
  flash memories.
\newblock {\em IEEE Transactions on Information Theory}, 59(11):7278--7287,
  2013.

\bibitem{ye2019some}
Z.~Ye, T.~Zhang, X.~Zhang, and G.~Ge.
\newblock Some new results on splitter sets.
\newblock {\em IEEE Transactions on Information Theory}, 66(5):2765--2776,
  2019.

\bibitem{zhang2016new}
T.~Zhang and G.~Ge.
\newblock New results on codes correcting single error of limited magnitude for
  flash memory.
\newblock {\em IEEE Transactions on Information Theory}, 62(8):4494--4500,
  2016.

\bibitem{zhang2017nonexistence}
T.~Zhang and G.~Ge.
\newblock On the nonexistence of perfect splitter sets.
\newblock {\em IEEE Transactions on Information Theory}, 64(10):6561--6566,
  2017.

\bibitem{zhang2023lattice}
T.~Zhang, Y.~Lian, and G.~Ge.
\newblock On lattice tilings of $\mathbb{Z}^n$ by limited magnitude error balls
  ${B}(n, 2, 1, 1)$.
\newblock {\em IEEE Transactions on Information Theory}, 69(11):7110--7121,
  2023.

\bibitem{zhang2017splitter}
T.~Zhang, X.~Zhang, and G.~Ge.
\newblock Splitter sets and $k$-radius sequences.
\newblock {\em IEEE Transactions on Information Theory}, 63(12):7633--7645,
  2017.




\end{thebibliography}
\end{document}